\documentclass{amsart}

\usepackage[dvips]{graphicx}

\usepackage{pgf,pgfarrows,pgfnodes,pgfshade}

\usepackage[T1]{fontenc}          
\usepackage{amssymb}

\theoremstyle{plain}
\newtheorem{mainthm}{Theorem}
\newtheorem{mainprop}[mainthm]{Proposition}
\newtheorem{maincorollary}[mainthm]{Corollary}

\newtheorem{theorem}{Theorem}[section]
\newtheorem{lemma}[theorem]{Lemma}

\theoremstyle{definition}
\newtheorem{definition}[theorem]{Definition}

\newtheorem{proposition}[theorem]{Proposition}

\theoremstyle{remark}
\newtheorem{remark}[theorem]{Remark}

\numberwithin{equation}{section}

\newcommand{\Z}{\mathbb{Z}}
\newcommand{\N}{\mathbb{N}}
\newcommand{\R}{\mathbb{R}}
\newcommand{\eps}{\varepsilon}

\newcommand{\dm}{\diff^{1}_{m}(M)}

\newcommand{\diff}{\operatorname{Diff}}

\begin{document}

\title[A generic condition for existence of symbolic extension]{A $C^1$ generic condition for existence of symbolic extensions of volume preserving diffeomorphisms}


\author{Thiago Catalan}





\maketitle

\begin{abstract}
We prove that a $C^1-$generic volume preserving diffeomorphism has a symbolic extension if and only if this diffeomorphism is partial hyperbolic. This result is obtained by means of good dichotomies. In particular, we prove Bonatti's conjecture in the volume preserving scenario. More precisely, in the complement of Anosov diffeomorphisms we have densely robust heterodimensional cycles. \end{abstract}

\section{Introduction and Statement of the Results}


A system $(X,f)$ has a {\it symbolic extension}, if there exist a subshift $(Y,\sigma)$, which is a closed, shift invariant subset of a full shift over a finite alphabet, and a surjective continuous map $\pi:Y\rightarrow M$ such that
$\pi\circ\sigma=f\circ\pi$.  In this case, $(Y,\sigma)$ is called an  {\it extension} of
$(X,f)$ and  $(X,f)$ a {\it factor} of $(Y,\sigma)$. 

One way to measure the complexity of a system $(X,f)$ could be by means of the topological entropy. Hence, if a system has a symbolic extension its complexity is bounded above by the complexity of a subshift. However, this system may contain additional information. The {\it symbolic extension entropy} of the system is the infimum of the topological entropy of all symbolic extensions of the system. And note that the topological entropy of a system is less than or equal to the symbolic extension entropy. The difference between these functions is called {\it residual entropy} and represents how entropy is hidden at finer and finer scales. 

A symbolic extension of $(X,f)$ is a {\it principal extension} if the map $\pi$ is such that
$h_{\nu}(\sigma)=h_{\pi_{*}\nu}(f)$ for every  $\sigma-$invariant measure
$\nu\in\mathcal{M}(\sigma | Y)$, where $h_{\nu}(\sigma)$ is the metric entropy of
$\sigma$ with respect to $\nu$. Note, the residual entropy is zero if the system has a principal extension.

Let $M$ be a Riemannian, connected and compact manifold. A diffeomorphism $f: M\rightarrow M$ is {\it asymptotically $h-$expansive} if
$$
\lim_{\eps\rightarrow 0}\sup_{x\in M}h(f|B_{\infty}(x,\eps))=0,
$$
where $B_{\infty}(x,\eps)=\{y\in M;\, d(f^j(x),f^j(y))<\eps \text{ for every } j\in \N\}$, and $h( . )$ denotes the topological entropy.  If there exists $\eps_0$ such that $\sup_{x\in M}h(f|B_{\infty}(x,\eps))=0$ for every $0<\eps<\eps_0$, then $f$ is called {\it  $h-$expansive.} 

Boyle, D. Fiebig,  U. Fiebig \cite{Boyle} using entropy structure showed that asymptotically $h-$expansive  diffeomorphisms have principal extension. Hence, if a diffeomorphism has no symbolic extension, it should have entropy hidden no matter how thin is the scale. In other words,  should there exist invariant subsets contained in balls with diameters arbitrary small having positive topological entropy. In some sense, this is in much the same way as the phenomena of coexistence of infinitely many horseshoes.

By a result of  Buzzi \cite{Buzzi}, every smooth diffeomorphism over a compact manifold is asymptotically $h-$ex\-pan\-si\-ve, then it has a principal extension. 
A conjecture of Downarowicz and Newhouse in \cite{DN} asserts that every $C^r$-diffeomorphism ($r\geq 2$) has a symbolic extension. 
This conjecture was solved for surface diffeomorphisms by Burguet \cite{Burguet}. Also, we would like to remark that recently Burguet and Fisher \cite{BF} extended this result to higher dimensions proving that every $C^2$ partially hyperbolic diffeomorphisms with a 2-dimensional center bundle has a symbolic extension.  

Hence, it seems natural to try to relate the existence of symbolic extensions to the differential structure of a diffeomorphism. 
For instance, using shadowing we can easily find a symbolic extension for Anosov diffeomorphisms. Diaz, Fisher, Pacífico and Vieitez \cite{diazfisher} showed that every $C^1$ partial hyperbolic diffeomorphism is $h-$expansive if we define partial hyperbolic diffeomorphisms as in \cite{CSY}. Hence it has a principal extension. See also \cite{df}. A diffeomorphism $f$ exhibits a {\it homoclinic tangency} if
there exists a hyperbolic periodic point $p$ of $f$ having a non transversal homoclinic point.  Thus, Liao, Viana and Yang proved that if a diffeomorphism is not approximated by one exhibiting homoclinic tangency, then it is also $h-$expansive.


On the other hand, we can consider a problem about the existence of diffeomorphisms that has no symbolic extensions. We can note that such diffeomorphisms have a rich dynamics, since they are not asymptotically $h-$expansive. Moreover, if the conjecture of Downarowicz and Newhouse is right such diffeomorphisms can not be $C^2$. 

In the symplectic scenario the author with Tahzibi \cite{CT} extended a result of  Downarowicz and Newhouse \cite{DN}, proving that $C^1-$generically either a symplectic diffeomorphism  is Anosov or has no symbolic extensions. That is, in the symplectic setting we have a large set of diffeomorphisms having no symbolic extensions.

The aim of this paper is to obtain similar results in the conservative case. 


We denote by $\dm$ the set of $C^1$ volume preserving diffeomorphisms over $M$. Here, as in \cite{CSY}, a $f-$invariant subset $\Lambda$ is {\it partial hyperbolic} if there exists a continuous $Df-$invariant splitting $T_{\Lambda}M=E^s\oplus E^c_1\oplus\ldots\oplus E^c_k\oplus E^u$, such that each center bundle $E_i^c$ is one-dimensional, and there exist constants $m\in\N$, $0<\lambda<1$ such that for 
every $x\in M$:
\begin{align*}
- &\|Df^m(v)\|\leq 1/2 \text{  for each unitary } v\in E^s \text{(one says } E^s \text{ is (uniformly) contracted.)},
\\
- &\| Df^{-m}(v)\|\leq 1/2 \text{  for each unitary } v\in E^u \text{(one says } E^u \text{ is (uniformly) expanded.)},
\\
- &\|Df^m_x(u)\|\leq 1/2 \|Df_x^m(v)\|, \text{ for each } x\in \Lambda, \text{ each } i=0,\ldots,k \text{ and each unitary}
\\
& \text{vectors } u\in E^s\oplus\ldots\oplus E^c_i, v\in E^c_{i+1}\oplus\ldots\oplus E^u \text{ in } T_xM.
\end{align*}
If all center bundles are trivial, then $\Lambda$ is called a {\it hyperbolic set}. 
We say that a volume preserving diffeomorphism $f:M\rightarrow M$ is {\it partial hyperbolic} if $M$ is a partial hyperbolic set. If $M$ is a hyperbolic set then we say that $f$ is an Anosov diffeomorphism.

The main result of this paper is the following:

 \begin{mainthm}\label{simboext}
There is a residual subset $\mathcal{R}\subset\diff^1_m(M)$ ($dim \, M\geq 3$) such that
if $f\in \mathcal{R}$ is a non partial hyperbolic diffeomorphism then $f$ has no symbolic extension.
\end{mainthm}

\begin{remark}
In dimension two, the previous theorem follows from  Downarowicz and Newhouse's result \cite{DN}.
\end{remark}

Now, Theorem \ref{simboext} and the result of Diaz, Fisher, Pacífico and Vieitez \cite{diazfisher} provide a generic intrinsic characterization of the existence of symbolic extensions in the volume preserving scenario.

\begin{mainthm}\label{symbequivalence}
There exists a residual subset  $\mathcal{R}\subset\diff^1_m(M)$ ($dim \, M\geq 3$) such that a diffeomorphism $f\in \mathcal{R}$ has a symbolic
extension if and only if it is partial hyperbolic. In particular, if $f\in\mathcal{R}$ has a symbolic extension then it has a principal extension. 
\end{mainthm}

A directly consequence of this result is the following.

\begin{maincorollary}
If a $C^1$ generic volume preserving diffemorphism $f$ is conjugated to a $C^{\infty}$ diffeomorphism, then $f$ is partial hyperbolic. 
\end{maincorollary}

In the articles \cite{DN} and \cite{CT} the main tool to obtain their results is the existence of an ``abundance'' of diffeomorphisms exhibiting homoclinic tangency in the complement of Anosov diffeomorphisms, since they are in the symplectic scenario.  
Hence, it is somewhat folklore the relation between robustness of homoclinic tangency and non existence of symbolic extensions.  

A diffeomorphism $f\in \diff^1(M)\, \big(\text{ resp. }  \diff_m^1(M)\big)$ exhibits a $C^1$ {\it robust
homoclinic tangency} if there exist a hyperbolic basic set $\Lambda$ of $f$ and a small
neighborhood $\mathcal{U}\subset \diff^1(M)\, \big(\text{ resp. } \diff^1_m(M)\big)$ of $f$ such that
$W^s(\Lambda(g))$ has a non transversal intersection with $W^u(\Lambda(g))$ for every
$g\in \mathcal{U}$. Where $\Lambda(g)$ is the continuation of $\Lambda$ for $g$.

\begin{mainprop} If $f\in \diff^1(M) \, \big(\text{ resp. } \diff^1_m(M)\big)$ exhibits a $C^1$
robust homoclinic tangency, then there exists a residual subset $\mathcal{R}$ in some
neighborhood of $f$ in $\diff^1(M)$ $ \big(\text{ resp. } \diff^1_m(M)\big)$, such that every $g\in \mathcal{R}$ has no symbolic extensions.
\label{mainprop}\end{mainprop}

As a consequence of this result we will obtain
Theorem \ref{simboext}. For that, we should investigate the relation between robutness of homoclinic tangency and partial hyperbolicity. More general, we should investigate the existence of good dichotomies.  

Recently, Crovisier, Sambarino and Yang \cite{CSY} proved that diffeomorphisms in $\diff^1(M)$ far from diffeomorphisms exhibiting homoclinic tangency, 
are approximated by partial hyperbolic diffeomorphisms. There they affirm that this result is also true in the conservative setting. Just for sake of completeness we will state it here and a sketch of the proof will appear inside the proof of Lemma \ref{lemaphtg}, see Remark \ref{csy}.

\begin{proposition}
Any diffeomorphism $f$ can be approximated in $\dm$ by diffeomorphisms which exhibit a homoclinic tangency or by partially hyperbolic diffeomorphisms.
\label{CSY}\end{proposition}

We define the {\it index} of a hyperbolic periodic point $p$ as the dimension of its stable manifold and we denote it by $ind\, p$.  A diffeomorphism $f$ exhibits a {\it heterodimensional cycle} if there exist hyperbolic periodic points $p$ and $q$ with different indices 
such that $W^s(p)\cap W^u(q)$ and $W^u(p)\cap W^s(q)$ are non empty intersections. 

One open problem about dichotomies is  Palis's conjecture, which says that densely in $\diff^r(M)$ ($r\geq 1$) either a diffeomorphism is hyperbolic, or exhibits a homoclinic tangency, or exhibits a heterodimensional cycle.
Palis's conjecture was proved for $C^1$ surface diffeomorphisms  by  Pujals and Sambarino \cite{PS}, and recently Crovisier and Pujals proved a remarkable result in this direction, by means of essential hyperbolicity,  for this result see \cite{CP}.  For symplectic and volume preserving diffeomorphisms,  there are also complete proofs for Palis's conjecture, see \cite{Newhouse2},   \cite{AC} and \cite{C}.
We would like to remark that in the volume preserving case what was proved, in fact, is that in the lack of hyperbolicity there are densely diffeomorphisms exhibiting heterodimensional cycles.  
Hence, our next theorem is a generalization of this result.  

A diffeomorphism $f\in \diff^1(M)\,$ \big(resp. $\diff_m^1(M)$\big) exhibits a  $C^1$ {\it robust
heterodimensional cycle} if there exist a hyperbolic basic set $\Lambda$ and a hyperbolic
periodic point $p$ of $f$, with $ind\, \Lambda\neq ind\, p$, such that $\Lambda(g)$ and $p(g)$
exhibits a heterodimensional cycle, i.e., $W^s(\Lambda(g))\cap W^u(p(g))$ and $W^u(\Lambda(g))\cap W^s(p(g))$ are non empty intersections, for every diffeomorphism $g$ in a
neighborhood of $f$ in $\diff^1(M)$ \big(resp. $\diff_m^1(M)$\big).


\begin{mainthm}\label{rcycle}
There is an open and dense subset $\mathcal{A}\subset\diff^1_m(M)$ ($dim \, M\geq 3$),
such that if $f\in \mathcal{A}$ is a non Anosov diffeomorphism, then $f$ exhibits a
$C^1$ robust heterodimensional cycle.
\end{mainthm}

It is worth to pointing out that the previous theorem is in fact  Bonatti's conjecture restrict to the volume preserving scenario. See  \cite{BD}. 


In order to prove Theorem \ref{rcycle}, we use the so known blenders. Bonatti and Diaz in \cite{BD} developed a way to obtain diffeomorphisms exhibiting robust homoclinic tangencies from blenders. In this paper, we develop their technics in the conservative case to prove the following result:

\begin{mainthm}\label{rtangency}
There is an open and dense subset $\mathcal{D}\subset \diff^1_m(M)$ ($dim \, M\geq 3$) such that if
$f\in\mathcal{D}$ is a non partial hyperbolic diffeomorphism then $f$ exhibits a robust homoclinic tangency.
\end{mainthm}

Note, Theorem \ref{simboext} is a directly consequence of the previous theorem and Proposition \ref{mainprop}.

This paper is organized as follows: in the second section we  recall some useful
perturbation results, and we show how to build a special kind of blender in the volume preserving setting. This special kind of blender is a blender horseshoe introduced in \cite{BD}. In this section we will prove Theorem \ref{rcycle}, too. In Section 3,  
we prove Theorem \ref{rtangency}, and 
finally, in Section 4, we prove Proposition \ref{mainprop} and Theorem \ref{simboext}.


\section{A blender horseshoe and some perturbation results}

Let $f$ be a $C^1$ diffeomorphism on $M$.  A hyperbolic transitive set $\Gamma$ of $f$  with 
$dim\, W^u(\Gamma,f)$ $=k\geq 2$
is   a {\it $ cu$-blender}  if there exist a $C^1$-neighborhood $\mathcal{U}$ of $f$ and a $C^1-$ open set $\mathcal{D}$ of embeddings of $(k-1)-$dimensional disks $D$ into $M$ such that, for every diffeomorphism $g\in \mathcal{U}$, every disk $D\in \mathcal{D}$ intersects the local stable manifold $W^s_{loc}(\Gamma(g))$, where $\Gamma(g)$ is the continuation of the hyperbolic set $\Gamma$ for $g$.  $\mathcal{D}$ is called the superposition region of the blender. Similarly, we can define a $cs$-blender with stable manifold replaced by unstable manifold.  The above definition was given in \cite{BD}. We would like to remark,  that for a $cu-$blender, these $(k-1)-$dimensional disks  are usually $uu-$disks. See Remark \ref{reblender}. 
We refer the reader to \cite{BD1} and \cite{BD} for more details about the geometry structure of this amazing set. 

The main result  of this section is the following.

\begin{proposition}\label{propblender}
If $f\in\dm$ has two hyperbolic periodic points $p_1$ and $p_2$ of different indices, say
$i$ and $i+j$, respectively, then for any neighborhood $\mathcal{U}\subset \dm$ of $f$
and any $i\leq k\leq i+j-1$ there exists a diffeomorphism $g\in \mathcal{U}$ having a $cu-$blender horseshoe  $\Gamma$ with index $k$.\end{proposition}

\begin{remark}
A $cu$-blender horseshoe is a special kind of a $cu-$blender set, which will be defined in the proof of Proposition \ref{propblender}.  
\end{remark}

\begin{remark}
A similar result still holds for $cs$-blenders. More precisely, if $f$ is under the same hypotheses of Proposition \ref{propblender}, then in any neighborhood of $f$ there exists a diffeomorphism $g$ having a $cs-$blender horseshoe of index $k$, for any $i+1\leq k\leq i+j$.
\label{csblender}\end{remark}

We also remark that F. Hertz, M. Hertz, Tahzibi and Ures have already showed in \cite{HHTU} the existence of blenders in the conservative scenario, after some perturbation, if the initial diffeomorphism is under the same hypotheses of Propostion \ref{propblender}. However, we observe that here we are interested in a special kind of blenders. Moreover, we emphasize that our methods to prove Proposition \ref{propblender} is different from their.


We will recall now some useful perturbation results. 

The first one is a Pasting lemma of Arbieto and Matheus \cite{pasting lema}.

\begin{theorem}[Pasting lemma] 
If $f$  is a  $C^2$ volume preserving diffeomorphism over $M$, and $x\in M$, then for every $\eps>0$ there exists a $C^1$ volume preserving diffeomorphism $g$ $\eps-C^1$ close to $f$ such that for a small neighborhood $U\supset V$ of $x$, $g|U^c=f$ and $g|V=Df(x)$ (in local coordinates). 
\label{pasting lemma}\end{theorem}

\begin{remark}
If $f$ is $C^{\infty}$ then $g$ could be taken $C^{\infty}$, too. 
\end{remark}

A directly consequence of pasting lemma is a conservative version of Franks lemma, see \cite{LLS}.

\begin{lemma}[Franks lemma]
\label{l.franks} Let $f\in \dm$ and $\mathcal{U}$ be a $C^1$ neighborhood of  $f$ in $\dm$. Then, there exists a smaller neighborhood
$\mathcal{U}_0\subset \mathcal{U}$ of $f$ and $\delta>0$ such that if $g\in \mathcal{U}_0(f)$, $S=\{x_1,\dots,x_m\}\subset M$ be any finite peace of orbit and
$\{L_i:T_{x_i}M\to T_{x_{i+1}}M\}_{i=1}^m$ conservative linear maps satisfying $\|L_i-Dg(x_i)\|\leq\delta$ for every
$i=1,\dots m$,  then for any small fixed neighborhood $V$ of $S$ there exist  $h\in\mathcal{U}(f)$ in the same class of differentiability  of $g$, such that $h=g$ in $V^c$, moreover $h(x_i)=g(x_i)$ and $Dh(x_i)=L_i$.
\end{lemma}

The next result is a connecting lemma of Hayashi \cite{H}. A conservative version was proved by Wen and Xia \cite{Xia-wen}. 

\begin{theorem}[$C^1$-connecting lemma]
Let $f \in \dm$ and $p_1,\, p_2$ hyperbolic periodic points of $f$, such that there exist sequences  $y_n\in M$ and positive integers 
$k_n$ such that:
\begin{itemize}
\item $y_n\rightarrow y \in W_{loc}^u(p_1, f))$, $y\neq p_1$; and
\item $f^{k_n}(y_n)\rightarrow x \in W_{loc}^s(p_2, f))$, $x\neq p_2$.
\end{itemize}
Then, there exists a $C^1$ volume preserving diffeomorphism  $g$   $C^1-$close to $f$ such that $W^u(p_1,g)$ and $W^s(p_2,g)$ have a non empty intersection close to  $y$. \label{connecting lema}\end{theorem}

The following technical result will be needed in the proof of Proposition \ref{propblender}.

\begin{lemma}\label{technicallemma}
If $f\in\dm$ has two hyperbolic periodic points $p_1$ and $p_2$ of different indices, say
$i$ and $i+j$, respectively, then for any neighborhood $\mathcal{U}\subset \dm$ of $f$,
$i\leq k\leq i+j-1$ and $\eps>0$ there exists a diffeomorphism $g\in \mathcal{U}$ with a hyperbolic
periodic point $p$, such that $p$ has index $k$, $Dg^{\tau(p,g)}$ has only real
eigenvalues with multiplicity one, say $\lambda_1<\ldots<\lambda_d$, and moreover
$|\lambda_{k+1}-1|<\eps$. Where $\tau(p,g)$ denotes the period of $p$ for $g$.\end{lemma}

This lemma follows by the same method as in Proposition 3.2 in \cite{AC}. However, provided this method will be useful later we will give a sketch of the proof.

Before we prove the above lemma, let us recall some definitions. 

Recall, two hyperbolic periodic points $p$ and $q$, having the same index are {\it homoclinically related} if there exist a transversal intersection between $W^s(p,f)$ and $W^u(q,f)$, and $W^u(p,f)$ and $W^s(q,f)$.  We denote by $H(p,f)$ the closure of the hyperbolic periodic points homoclinically related to $p$, which is called the {\it
homoclinic class} of $p$. Similarly, we can define when hyperbolic periodic points and hyperbolic sets are homoclinically related. 

A continuous $Df-$invariant splitting $T_{\Lambda}M= E_1\oplus\ldots\oplus E_k$ for a $f-$invariant subset $\Lambda$ is {\it dominated} if the third condition in the partial hyperbolic definition is satisfied. 

For abbreviation, sometimes we use expressions like "after a perturbation", or "there exists a diffeomorphism $C^1-$close", which means that these perturbations could be done so small as we want. 

\vspace{0,2cm}
{\it Proof of Lemma \ref{technicallemma}:}

By a result of Xia \cite{Xia},
a generic volume preserving diffeomorphism has all homoclinic classes non trivial. Thus, after a perturbation, we can assume $H(p_1,f)$ and $H(p_2,f)$ are non trivial.  Now, by results of Bonatti, Diaz and Pujals \cite{BDP}, and Franks lemma  we can perturb $f$ to $f_1$ in order to obtain $\tilde{p}_1$ and
$\tilde{p}_2$ hyperbolic periodic points homoclinically related to $p_1(f_1)$ and $p_2(f_1)$, respectively, such that $Df_1^{\tau(\tilde{p}_1,f_1)}(\tilde{p}_1)$ and
$Df_1^{\tau(\tilde{p_2},f_1)}(\tilde{p_2})$ have only real eigenvalues with multiplicity one. 

To simplify the notation we replace $\tilde{p}_1$ and $\tilde{p}_2$ by $p_1$ and $p_2$, respectively. 
And moreover, we will continue to write  $p_1$ and $p_2$ for their continuations.  

In the sequence we perturb $f_1$ in order to find a diffeomorphism exhibiting a heterodimensional cycle between $p_1$ and $p_2$. For that, we use a result of
Bonatti and Crovisier \cite{BONATTICROVISIER}:

\begin{proposition}[Bonatti and Crovisier] 
There exists a residual subset $\mathcal{R}$ of $\dm$ such that if $g\in \mathcal{R}$ then there exists a hyperbolic periodic point $p$ of $g$ such that $M=H(p,g)$. In particular, $g$ is transitive. \label{BC}\end{proposition}

After a perturbation, we can assume  $f_1\in \mathcal{R}$, i.e, $f_1$ is transitive. Then, using connecting lemma
we can perturb $f_1$ to $f_2$ such that there is a non transversal intersection between $W^u(p_1,f_2)$ and $W^s(p_2,f_2)$.  Since this intersection is robust, we can repeat the above process and perturb $f_2$ to $f_3$  such that $W^s(p_1,f_3)$ and $W^u(p_2,f_3)$ have also a non empty intersection,  which implies that $f_3$ exhibits a heterodimensional cycle between $p_1$ and $p_2$. Moreover, $f_3$
can be taken such that  $Df_3^{\tau(p_1,f_3)}(p_1)$ and $Df_3^{\tau(p_2,f_3)}(p_2)$ have only real eigenvalues with multiplicity one.

Let $x\in W^s(p_1,f_3)\cap W^u(p_2,f_3)$ and $y\in W^u(p_1,f_3)\cap W^s(p_2,f_3)$ be two heteroclinic
points of the cycle. Recall $ind\, p_1=i$ and $ind\, p_2=i+j$. Without loss of generality we can assume $y$ is a transversal heteroclinic point, and $x$ is a quasi-transversal heteroclinic point, i.e., $T_xW^s(p_1,f_3)\cap T_x W^u(p_2,f_3)=\{0\}$.  By the regularization result of Ávila \cite{Avila} (which says we can suppose $f_3$ be $C^{\infty}$) and pasting lemma, we can linearize the diffeomorphism in a
small neighborhood $U_{p_1}$ and $U_{p_2}$ of  $p_1$ and $p_2$, respectively.  More precisely,  we can assume  $f_3$ is equal to $Df_3(p_1)$ and $Df_3(p_2)$ (in local coordinates) in the neighborhoods $U_{p_1}$ and $U_{p_2}$, respectively. 

For simplicity of notation, in the reminder of this proof we assume that $p_1$ and $p_2$ are fixed points, and we will look at $U_{p_1}$ and $U_{p_2}$ in local coordinates. Since $Df_3(p_1)$ and $Df_3(p_2)$ have only real eigenvalues with multiplicity one, we can find a decomposition of $\R^d$ by eigenspaces of $Df_3(p_1)$ (resp. $Df_3(p_2)$), which we denote by $E_{1,p_1}\oplus\ldots\oplus E_{d,p_1}$ (resp. $E_{p_2}\oplus\ldots\oplus E_{d,p_2}$). 
We set $\lambda_k$ (resp. $\sigma_k$), $k=1,\ldots, d$, the eigenvalue of $Df_3(p_1)|E_{k,p_1}$ (resp. $Df_3(p_2)|E_{k,p_2}$). We can also suppose the eigenvalues are in an increase order. 

In order to be more precise, we will make the following assumptions, we consider $E_{i,p_1}(.)$ the extension of the direction $E_{i,p_1}$ in the neighborhood $U_{p_1}$, the same for $E_{i,p_2}(.)$. We remark these decompositions are all dominated splittings, indeed.




 $ind\, p_1=i$, it follows that the stable and unstable directions of $p_1$ are $E^{s}_{p_1}=E_{1,p_1}(p_1)\oplus \ldots \oplus E_{i,p_1}(p_1)$ and $E^{u}_{p_1}=E_{i+1,p_1}(p_1)\oplus
\ldots \oplus E_{d,p_1}(p_1)$,  respectively. Similarly,  the stable and unstable directions of $p_2$ are $E^{s}_{p_2}=E_{1,p_2}(p_2)\oplus \ldots \oplus E_{i+j,p_2}(p_2)$ and
$E^{u}_{p_2}=E_{i+j+1,p_2}\oplus \ldots \oplus E_{d,p_2}$, respectively, since $ind\, p_2=i+j$.  Moreover, by the choice of $f_3$, 
$W^s_{loc}(p_1,f_3)=E^s_{p_1}\cap U_{p_1}$, $W^u_{loc}(p_1,f_3)=E^u_{p_1}\cap U_{p_1}$, $W^s_{loc}(p_2,f_3)=E^s_{p_2}\cap U_{p_2}$ and $W^u_{loc}(p_2,f_3)=E^u_{p_2}\cap U_{p_2}$.


\vspace{0,3cm}
{\bf Claim:} {\it There is a diffeomorphism $f_4$ $C^1-$close to $f_3$ such that the $f_3$-invariant subset $\Lambda=O(x)\cup O(y)\cup\{p_1,p_2\}$ still is $f_4-$invariant and moreover has a dominated splitting by one dimensional sub-bundles. }

\vspace{0,2cm}

We define $E(y):=T_y (W^u(p_1,f_3)\cap W^s(p_2,f_3))$. 
Since $y$ belongs to  unstable manifold of $p_1$, and $f_3|U_{p_1}=Df_3(p_1)$,  if $n$ is large enough, it follows that  $Df_3^{-n}(y)(E(y))$ is in $E_{i+1,p_1}(f_3^{-n}(y))\oplus\ldots \oplus E_{d,p_1}(f_3^{-n}(y))$. 
Moreover, 
by transversality 
we can assume that $Df_3^{-n}(y)(E(y))\cap
E_{i+j+1,p_1}(f_3^{-n}(y))\oplus \ldots\oplus E_{d,p_1}(f_3^{-n}(y))=\{0\}$. Provided we have a dominated splitting in $U_{p_1}$, $Df_3^{-n}(y)(E(y))$ converges to
$E_{i+1,p_1}(p_1)\oplus\ldots\oplus E_{i+j,p_1}(p_1)$ when $n\to\infty$. Then, choosing $n$ large enough and using Franks lemma, after a perturbation we can assume $f_3$
such that $Df_3^{-n}(y)(E(y))=E_{i+1,p_1}(f_3^{-n}(y))\oplus\ldots\oplus E_{i+j,p_1}(f_3^{-n}(y))$. Note, the perturbation necessary here is local, and moreover keeps unchanged the orbit of $y$.


We now apply this argument again, considering the future orbit of $y$, to obtain a perturbation of $f_3$ such that we have also
$Df_3^{n}(y)(E(y))=E_{i+1,p_2}(f_3^n(y))\oplus\ldots \oplus E_{i+j,p_2}(f_3^n(y))$. This perturbation of $f_3$ which we continue denoting by the same letter has a $Df_3$-invariant subbundle on $O(y)\cup \{p,q\}$ which we will denote by $E$, for convenience.


The $\lambda-$lemma says that $Df_3^{m}(f_3^{-n}(y))(E_{i+j+1,p_1}(f_3^{-n}(y))\oplus \ldots \oplus E_{d,p_1}(f_3^{-n}(y))$ converges to
$E_{i+j+1,p_2}(p_2)\oplus\ldots \oplus E_{d,p_2}(p_2)$  if  $m\to\infty$. Then by the same argument again we can perturb $f_3$ such that $Df_3^{m}(f_3^{-n}(y))(E_{i+j+1,p_1}(f_3^{-n}(y))\oplus \ldots \oplus E_{d,p_1}(f_3^{-n}(y))=$  $ E_{i+j+1,p_2}(f_3^{m-n}(y))\oplus \ldots \oplus E_{d,p_2}(f_3^{m-n}(y))$, and the sub-bundle $E$ is still $Df-$invariant. Replacing $m$ and $n$ by large positive integers if necessary, and applying once more the argument, $f_3$ could be assumed such that
$Df_3^{-m}(f_3^{m-n}(y))(E_{1,p_2}(f_3^{m-n}(y))\oplus \ldots \oplus E_{i,p_2}(f_3^{m-n}(y))=E_{1,p_1}(f^{-n}(y))\oplus \ldots \oplus E_{i,p_1}(f^{-n}(y))$. 

Therefore, $f_3$ is such that there exists a $Df_3-$invariant splitting over $O(y)\cup\{p_1,p_2\}$.

Moreover, if we repeat this process finitely many times inside each invariant sub-bundle, $f_3$ could be assumed such that
$$
Df_3^{2n}(f_3^{-n}(y)(E_{k,p_1}(f_3^{-n}(y)))=E_{k,p_2}(f_3^{n}(y)), \, k=0,\ldots, d; \text{ for }\, n \text{ large enough}.
$$

Finally, applying the above arguments, with $y$ replaced by $x$,  $f_3$ can also be chosen such that 
$$
Df_3^{2n}(f_3^{-n}(x))(E_{k,p_2}(f_3^{-n}(x))=E_{k,p_1}(f_3^{n}(x)), \, k=0,\ldots, d; \text{ for }\, n \text{ large enough},
$$
which finishes the proof of the claim, since this $Df-$invariant splitting is natural dominated.
 

\vspace{0,3cm}
We fix now an arbitrary $i\leq k\leq i+j-1$, and we consider the diffeomorphism $f_4$ given by the previous Claim. 
Using the heteroclinic points $x$ and $y$, we can perform a perturbation of $f_4$  to obtain a periodic orbit  in a small neighborhood of  $\Lambda$, with  arbitrary large period. In fact, this could be done such that this periodic orbit has as many points as we want in the neighborhoods $U_{p_1}$ and $U_{p_2}$, being fixed the number of points outside these neighborhoods. Hence, by  continuity of the dominated splitting over $\Lambda$, and since $\|Df_4| E_{k+1}(p_1)\|>1$ and $\|Df_4| E_{k+1}(p_2)\|<1$, it follows there exists a diffeomorphism $f_5$ $C^1-$close to $f_4$ having a hyperbolic periodic point $p$ in a neighborhood of $\Lambda(f_5)$ with index $k$ and such that $Df_5^{\tau(p,f_5)}(p)$ has only real eigenvalues with multiplicity one. Moreover, this could be done such that $\|Df_5^{\tau(p,f_5)}| E_{k+1}(p)\|$ is so close to one as we want. For details we refer the reader to \cite{AC}.


$\hfill\square$

{ \it Proof  Proposition \ref{propblender} :} 


We fix an arbitrary $i\leq k\leq i+j-1$.  By Lemma \ref{technicallemma}, after a perturbation,  we can assume there exists a hyperbolic periodic point $p$ of $f$ such that $p$ has index $k$, $Df^{\tau(p,f)}(p)$ has only real eigenvalues with multiplicity one, say $\lambda_1<\ldots<\lambda_d$, and moreover $\lambda_{k+1}$ is so close to one as
we want.

If  $E_{\lambda_t}$ is the corresponding eigenspace to $\lambda_t$, then we have on $p$ a natural partially hyperbolic splitting $T_pM=E^{s}\oplus E^{cu}\oplus
E^{uu}$, where $E^s=\cup_{1\leq t\leq k} E_{\lambda_t}$ is the stable direction of dimension
$k$, and the unstable direction is divided in two subspaces,
$E^{cu}=E_{\lambda_{k+1}}$ (the center unstable direction), 
and $E^{uu}=\cup_{t> k+1} E_{\lambda_{t}}$ the strong unstable direction. By Hirsch,  Pugh and  Shub \cite{HPS},  the strong directions are integrable, which means here the existence of $W^{uu}(p,f)$, the { \it strong unstable manifold} of $p$, which varies continuously with respect to the diffeomorphism.

As in the proof of Lema \ref{technicallemma}, 
we can perturb $f$ to a $C^{\infty}$ diffeomorphism $f_1$ to obtain a intersection  
between the stable and strong unstable manifolds of $p(f_1)$,  $W^{s}(p(f_1),f_1)\cap W^{uu}(p(f_1),f_1)\neq \emptyset$, and moreover such that $f_1^{\tau(p_1(f_1),f_1)}=$\newline $Df^{\tau(p_1(f_1),f_1)}_1(p(f_1))$ (in local coordinates) in
a neighborhood $U$ of $p(f_1)$. 

By abuse of notation, we write just $p$ instead of $p(f_1)$.
Also, since $\|Df^{\tau(p,f)}(p)|E^{cu}\|$ is so close to one as we want, after another perturbation, we can suppose
$|\tilde{\lambda}_{c}|=\|Df_1^{\tau(p,f_1)}(p)|E^{cu}\|=1$. 

From now on, we look at $U$ in local coordinates. Then,  in $U$ the local stable and
strong unstable manifolds of $p$ coincide with their directions, i.e.,
$W^s_{loc}(p,f_1)=E^s(p,f_1)\cap U$ and $W^{uu}_{loc}(p,f_1)=E^{uu}(p,f_1)\cap U$.

Let $x\in W^{s}(p,f_1)\cap W^{uu}(p,f_1)$ be a strong  homoclinic point of $p$, which we can also assume be a
quasi-transversal strong homoclinic point, i.e.,
$dim(T_xW^s(p,f_1)+T_xW^{uu}(p,f_1))=d-1$ since $dim\, E^{cu}(p)=1$. Hence, there exist positive integers
$n$ and $m$ such that $f_1^n(x)=(x^s_0,0,0)$, $f_1^{-m}(x)=(0,0,x^u_0)\in U$.
Here, we are considering the natural extension to $U$ of the partial hyperbolic splitting $T_pM=E^s\oplus E^{cu}\oplus E^{uu}$. 
Also, without loss of generality we can suppose this decomposition orthogonal. 

By the same method as in the Claim in the proof of Lemma \ref{technicallemma}, we can find a diffeomorphism $f_2$ $C^1-$close to $f_1$, such that shrinking $U$ if necessary $f_2$ satisfies the
following conditions:

\begin{itemize}
\item[1-] $f_2^{\tau(p,f_2)}=Df_2^{\tau(p,f_2)}(p)=Df_1^{\tau(p,f_1)}(p)$  in $U$, keeping invariant the directions $E^j\cap U$, $j=s,\,cu,\,
uu$; 

\item[2-] $x$ is still a strong homoclinic point of $p$, and moreover $$Df_2^{mn}(f_2^{-m}(x))(E^j(f_2^{-m}(x)))=E^j(f^n(x)), \; j=s,\, cu,
\, uu.$$
\end{itemize}

 $f_2$ is obtained through several finitely many perturbations of $f_1$ using Franks lemma,
$f_2$ is in the same class of differentiability of $f_1$, which implies 
$f_2$ is $C^{\infty}$.  Hence, we can use pasting lemma in order to linearize $f_2$ 
in a segment of the orbit of $x$.  More precisely, we can choose neighborhoods $U_m,\, U_n\subset U$ of
$f_3^{-m}(x)$ and $f_3^{n}(x)$, respectively, and perturb $f_2$ to $f_3$  such that  $f_3^{mn}(E^j(y)\cap U_m)=E^j(f^{mn}(y))\cap U_n$, for every $y\in U_m$ and $j=s,\, cu,\, uu$. Using the fact that $f_2^{\tau(p,f_2)}$ is linear in $U$ and
$\tilde{\lambda}_c=1$, replacing $m$ and $n$ with larger ones, and after once more perturbation, we can suppose $f_3$ satisfying

\begin{itemize}
\item[3-]$f_3^{mn}:U_m\rightarrow U_{n}$ is an affine map.  More precisely,
$$f_3^{mn}(x^s, x^c,x^u)=(x^s_0+A_s(x^s),
\lambda_c x^c, A_u(x^u-x^u_0)),$$ where $A_s$ is a linear contraction, $A_u$ a linear
expansion and $1<\lambda_c<1+\epsilon$, for some small $\eps>0$.

\item[4-] $E^s\oplus E^{uu}$ is invariant for both maps $f_3^{\tau(p,f_3)}|U$ and
$f_3^{mn}|U_m$.
\end{itemize}


Hence, if $D\subset (E^s\oplus E^{uu})\cap U$
is a small enough rectangle containing $p$ and $f_3^n(x)$ in its interior, then $f_3^{l\tau(p,f_3)+mn}(D)\cap D$ has two non empty disjoint connected components for some $l$ large
enough. 
One of them containing $p$
and another one $f_3^n(x)$, which we denote by $\mathcal{A}$ and $\mathcal{B}$,
respectively.

For simplicity of notation, we set $\tilde{F}=f_3^{l\tau(p,f_3)+mn}|D$,
$\mathbb{A}=\tilde{F}^{-1}(\mathcal{A})$ and $\mathbb{B}=\tilde{F}^{-1}(\mathcal{B})$.
Note, 
$\tilde{F}$ is a linear map on
$\mathbb{A}\cup \mathbb{B}$, and the stable and strong unstable directions are $\tilde{F}-$invariant. Moreover, taking $l$ larger if necessary  $\tilde{F}|E^s$ and
$\tilde{F}^{-1}|E^{uu}$ are contractions, for every point in
$\mathbb{A}\cup\mathbb{B}$ and $\mathcal{A}\cup \mathcal{B}$, respectively. Hence, the
maximal invariant set in $D$ for $\tilde{F}$,
$$
\Sigma=\bigcap_{j\in \Z} \tilde{F}^j(D)
$$
is a hyperbolic set conjugated to  the full shift of two symbols. We denote by
$q\in \mathbb{B}$ the other fixed point of $\tilde{F}$. Note, $E^s\oplus E^{uu}$ is the hyperbolic splitting over $\Sigma$.

Fixing any arbitrary small $\delta>0$, we set $R=D\times [ -\delta,\delta]\subset U$, and
replace $\mathbb{A}$ and $\mathbb{B}$ by $\mathbb{A}\times[-\delta,\delta]$ and
$\mathbb{B}\times[-\delta,\delta]$, respectively.  Taking $\delta$ smaller, 
$F:=f_3^{l\tau(p,f_3)+mn}|\mathbb{A}\cup \mathbb{B}$ is then well defined. Moreover, taking
the center coordinate as the last one, we have 
$$
F(x^s,x^u,x^c)=(\tilde{F}(x^s,x^u),
\lambda_c x^c).
$$
Since $\lambda_c>1$, it follows that  $\Lambda_0=\Sigma\times 0$ is the maximal
$F-$invariant set in $R$. Also, provided  $E^s$, $E^{cu}$ and $E^{uu}$ are $F$-invariant, we have a natural partial hyperbolic splitting on $\Lambda_0$. In particular, 
$\Lambda_0$ is a hyperbolic set with index $k$ since $\|F|E^{cu}\|>1$. 

After a coordinate change, we can suppose $R=[-1,1]^s\times [-1,1]\times
[-1,1]^u$ in local coordinates, and $p=(0,0,0)$ in this chart.  

For every $t>0$ small enough, using pasting lemma we can find a perturbation $h_t$ of the identity map such that  $h_t(x^s,x^c,x^u)=(x^s,x^c-t,x^u)$ for every point
in $U_n$ and $h_t=Id$ outside a small neighborhood of $U_n$.  We
define $ f_t=h_t\circ f_3, $ which is $C^1-$ close to $f_3$.

Shrinking $U_n$ if necessary, the above perturbation $f_t$ in terms of $F_t$ is
the following
$$
\begin{array}{llll}
1- & F_t=F, & \text{if } x\in
\mathbb{A}
\\
2- & F_t=F+(0,-t,0), & \text{ if } x\in
\mathbb{B}.
\end{array}
$$

Provided $t$ is small, the maximal $F_t-$invariant set $\Lambda_t$  in $R$  is the continuation of the hyperbolic set $\Lambda_0$ of $F$, hence $\Lambda_t$ is also
hyperbolic. Moreover, note $E^s\oplus E^{cu}\oplus E^{uu}$ is still the hyperbolic
splitting on $\Lambda_t$, and $p$ is still a hyperbolic fixed  point of $F_t$. We denote
by $q_t$ the continuation of the hyperbolic fixed point $q$ of $F$. 

This set $\Lambda_t$ is defined as a {\it $cu-$blender horseshoe}. 

In the sequence, we will describe some properties of $\Lambda_t$ which
characterize, in fact, a cu-blender horseshoe.

For  $\alpha\in (0,1)$ we denote by $C^s_{\alpha}$    
and $C^{uu}_{\alpha}$ the following cone-fields in $R$:
$$
\begin{array}{lll}
C^s_{\alpha}(x)&=\{v=(v^s,v^{c},
v^{u})\in E^s\oplus E^{cu}\oplus
E^{uu}=T_xM; & \|v^c+v^u\|\leq
\alpha\|v^s\|\}
\\
C^{uu}_{\alpha}(x)&=\{v=(v^s,v^{c},
v^{u})\in E^s\oplus E^{cu}\oplus
E^{uu}=T_xM; & \|v^s+v^c\|\leq
\alpha\|v^u\|\}.
\end{array}
$$
We say that a disk $\Delta$ of dimension $s$ contained in $R$ is a {\it $s-$disk}
if
\begin{itemize}
\item it is tangent to
$C^s_{\alpha}$, i.e., $T_x\Delta\subset
C^s_{\alpha}(x)$ for all $x\in \Delta$,
and

\item its boundary $\partial \Delta$ is
contained in $\{-1,1\}^s\times
[-1,1]\times [-1,1]^u$.
\end{itemize}
On the other hand, a disk $\Upsilon$ of
dimension $u$ is a {\it $uu$-dis}k if
\begin{itemize}
\item it is tangent to
$C^{uu}_{\alpha}$, i.e.,
$T_x\Upsilon\subset C^{uu}_{\alpha}(x)$
for all $x\in \Upsilon$, and

\item its boundary $\partial \Upsilon$
is contained in $[-1,1]\times
[-1,1]\times \{-1,1\}^u$
\end{itemize}

\begin{remark} As expected, $W_{loc}^s(p,F_t)=[-1,1]^s\times
\{0\}\times\{0\}$ and $W_{loc}^s(q_t,F_t)=[-1,1]^s\times \{t/\lambda_c(\lambda_c-1)\}\times
\{x_0^u(q)\}$  are natural $s-$disks, while $W^{uu}_{loc}(p,F_t)=\{0\}\times
\{0\}\times[-1,1]^u$ and $W^{uu}_{loc}(q_t,F_t)=\{x^s_0(q_t)\}\times
\{t/\lambda_c(\lambda_c-1)\}\times [-1,1] $ are natural $uu-$disks.
\end{remark}

Note, there are two different homotopy classes of $uu-$disks disjoint from
$W^s(p,F_t)$. We say that an $uu-$disk is at the right of $p$ if it belongs to the
same homotopy class of $W^{uu}(q_t,F_t)$, and at the left otherwise. Similarly, we say that  a
$uu-$disk is at the left of $q_t$ if it belongs to the  same homotopy class of
$W^{uu}(p,F_t)$, and at the right otherwise.

By this convention, if $\mathcal{D}$ is an $uu-$disk, then one of the following is true:

\begin{itemize}
\item[-] $\mathcal{D}$ is at the left of $p$;

\item[-] $\mathcal{D}\cap W^s(p,F_t)\neq
\emptyset$;

\item[-] $\mathcal{D}$ is at the right of $q_t$;

\item[-] $\mathcal{D}\cap W^s(q_t,F_t)\neq
\emptyset$;

\item[-] $\mathcal{D}$ is at the right of $p$ and
at the left of $q_t$. In this case we say
that the $uu-$disk is in between of $p$
and $q_t$.
\end{itemize}

We fix now a very small $\alpha\in (0,1)$ in the definition of the $uu$-disks, such that a
$uu-$disk $\mathcal{D}$ is $C^1-$close to $E^{uu}$. If we define
$F_{\mathbb{A}}(\mathcal{D})=F_t(\mathbb{A}\cap \mathcal{D})$ and $F_{\mathbb{B}}(\mathcal{D})=F_t(\mathbb{B}\cap \mathcal{D})$, then
the following is true:

\begin{itemize}

\item[1-] If $\mathcal{D}$ is at the right of $p$
(resp. $q_t$) then $F_{\mathbb{A}}(\mathcal{D})$
(resp. $F_{\mathbb{B}}(\mathcal{D})$) also is.

\item[2-] If $\mathcal{D}$ is at the left of $p$
(resp. $q_t$) then $F_{\mathbb{A}}(\mathcal{D})$
(resp. $F_{\mathbb{B}}(\mathcal{D})$) also is.

\item[3-] If $\mathcal{D}$ is at the left of $p$
 then $F_{\mathbb{B}}(\mathcal{D})$ also is.

\item[4-] If $\mathcal{D}$ is at the right of $q_t$
 then $F_{\mathbb{A}}(\mathcal{D})$ also is.

\item[5-] If $\mathcal{D}$ is in between of  $p$
and $q_t$, then either
$F_{\mathbb{A}}(\mathcal{D})$ or
$F_{\mathbb{B}}(\mathcal{D})$ is in between of
$p$ and $q_t$.
\end{itemize}

\begin{remark}
The above properties are robust. More precisely, if  $g$ is a diffeomorphism  $C^1-$close to $f_t$, and if we denote by $\Lambda_g$ the
continuation of the the hyperbolic periodic set $\Lambda_t$ of $f_t$, then
$g^{l\tau(p(g),g)+mn}|\Lambda_g=G|\Lambda_g$ has the same properties of $F_t|\Lambda_t$. Therefore, a $cu$-blender horseshoe set is robust. 
\end{remark}

\begin{remark}
Using iterated functions is possible to verify that every $uu$-disk in between of $W^s(p(g),g)$ and $W^s(q_t(g), g)$ intersects $W_{loc}^s(\Lambda_g,g)$. See Bonatti and Diaz \cite{BD1}.  In particular, the blender horseshoe $\Lambda_g$ is in fact a $cu$-blender, where the $uu$-disks in between of $W^s(p(g),g)$ and $W^s(q_t(g),g)$ define its superposition region.  
\label{reblender}\end{remark}


To see more properties about a blender horseshoe set we refer the reader to \cite{BD}. 

$\hfill\square$


We prove now Theorem \ref{rcycle}.

\vspace{0,3cm}
{\it Proof of Theorem \ref{rcycle}:}
\vspace{0,1cm}

Let $f$ be a non Anosov volume preserving diffeomorphism. By Theorem 1.1 in \cite{AC},  there exists a diffeomorphism $f_1\in \dm$ $C^1-$close to $f$ having a non-hyperbolic periodic point $p$. After a bifurcation of $p$ we can assume that $f_1$ has two hyperbolic periodic points of different indices, say $p_1$ and $p_2$, with  $ind\, p_1=i$ and $ind\, p_2=i+j$, $i,\, j>0$. 

By Proposition \ref{propblender} we can find a volume preserving diffeomorphism $f_2$ $C^1-$close to $f_1$ such that $f_2$ has a blender horseshoe $\Lambda$ with index $i+j-1$. We replace now $p_1$ by one of the two reference saddles of $\Lambda$. 

As in the proof of Proposition \ref{propblender},  we can perturb $f_2$ to $f_3$ to obtain a heterodimensional cycle between $p_1$ and $p_2$. Let $z$ denote a point of non transversal intersection between $W^s(p_1,f_3)$ and $W^{u}(p_2,f_3)$, which we can assume to be a quasi transversal intersection. 
Provided the partial hyperbolic structure in the superposition region $\mathbb{C}$ of the blender,  replacing $z$ by a positive iterated,  the connected disc in $W^u(p_2,f_3)\cap \mathbb{C}$ containing $z$ is in fact a $uu-$disk which is in between of the two reference saddles of $\Lambda$, as defined in the proof of Proposition \ref{propblender}. Note, this could be done such that $W^u(p_1,f_3)\cap W^s(p_2,f_3)$ has a transversal intersection. 

Therefore,  by properties of blenders, Remark \ref{reblender}, and continuity of the unstable manifold of $p_2$, every volume preserving diffeomorphism $g$ in a small neighborhood of $f_3$ has a heterodimensional cycle between $p_2(g)$ and $\Lambda(g)$.

$\hfill\square$

\section{Robustness of homoclinic tangency}

In this section we prove Theorem \ref{rtangency}. For that, it will be necessary to introduce folded submanifolds, introduced by Bonatti and Diaz in \cite{BD}.

\begin{definition}
Let $f$ be a diffeomorphism on $M$ having a blender-horseshoe set $\Lambda$ of index $u+1$ with reference cube $\mathbb{C}$, reference saddles $p$ and $q$, and $N\subset M$ be a submanifold of dimension $u+1$. We say that $N$ is {\it folded with respect to} $\Lambda$ if the interior of $N$ contains a sub-manifold $\mathcal{S}\subset \mathbb{C}\cap N$ of dimension $u+1$,  satisfying the following properties:

\begin{itemize}
\item There are $0<\alpha'<\alpha$ and a family $(S_t)_{t\in [0,1]}$ of  disks tangent to the cone field $C^{uu}_{\alpha'}$,  depending continuously on $t$, such that  $\mathcal{S}=\cup_{t\in [0,1]} \mathcal{S}_t$. Here,   $\alpha$ comes from the definition of a blender horseshoe, in particular $S_t$ is a $uu-$disk;

\item $\mathcal{S}_0\cap W^s_{loc}(A)$ and $\mathcal{S}_1\cap W^s_{loc}(A)$ are non empty transverse intersection points between $N$ and $W^s_{loc}(A)$, where  $A\in \{p,q\}$.

\item for every $t\in (0,1)$, the $uu-$disk $\mathcal{S}_t$ is in  between of $W^s_{loc}(p)$ and $W^s_{loc}(q)$.
\end{itemize}
To emphasize the reference saddle $A$ of the blender  we have considered, we say a submanifold $N$ is folded with respect to $(\Lambda,\, A)$.
\end{definition}


\begin{theorem}[Theorem 2, pg 18, \cite{BD}]\label{foldedmanifold}
Let $f$ be a $C^r$ ($r\geq 1$) diffeomorphism over $M$, and  $N\subset M$ be a folded submanifold with respect to a blender-horseshoe $\Lambda$ of $f$. Then $N$ and $W^s_{loc}(\Lambda)$ have a non empty $C^r-$robust non transversal intersection.
\end{theorem}

To prove Theorem \ref{rtangency} the following results will also be needed. 

First, Wen has proved in \cite{W} a dichotomy between diffeomorphisms having a dominated splitting and diffeomorphisms exhibiting a homoclinic tangency in the space of $C^1$ diffeomorphisms. Using the pasting lemma, Liang, Liu and Sun \cite{LLS} proved this dichotomy in the volume preserving scenario.  

\begin{proposition}\label{TvsD}[Theorem 1.3 in \cite{LLS}]
Let $f\in \dm$, and  $p$ be a hyperbolic periodic point of $f$. Then, we have the following dichotomy:
\begin{itemize}
\item [1-] Either the homoclinic class of $p$, $H(p,f)$, has a dominated splitting  $TM=E\oplus F$, with $dim\, E=ind(p)$, or

\item[2-] there exists a diffeomorphism $g$ $C^1-$close to $f$, exhibiting a homoclinic tangency for $p(g)$. 
\end{itemize}
\end{proposition}

The following result is a conservative version of  Theorem 1 in \cite{ABCDW}.  It may be proved using the same arguments as in Lemma \ref{technicallemma}, see \cite{AC} for details. 

\begin{proposition}
There is a residual subset $\mathcal{R}\in \dm$ of diffeomorphisms $f$ such that, for every $f\in \mathcal{R}$ containing hyperbolic periodic points of indices $i$ and $j$ contains hyperbolic periodic points of index $k$ for all $i\leq k\leq j$.
\label{intind}\end{proposition}


{\it Proof of Theorem \ref{rtangency}:}
\vspace{0,1cm}

Let $f$ be a volume preserving diffeomorphism which is not approximated by a partial hyperbolic diffeomorphism in $\dm$. In particular, $f$ is a non Anosov diffeomorphism, and then after a perturbation if necessary as in the proof of Theorem \ref{rcycle} we can assume $f$ has hyperbolic periodic points of different indices.

 We set $i$ and $j$ the smallest and largest positive integers, respectively, such that every hyperbolic periodic point $p$ of $f$ has $i\leq ind\, p\leq j$. 
  we are in the volume preserving scenario, after a perturbation, we can suppose this is still true for diffeomorphisms in a small neighborhood $\mathcal{U}$ of $f$. More precisely, if $g\in \mathcal{U}$ and $p$ is a hyperbolic periodic point of $g$ then $i\leq ind\, p\leq j$.

By Proposition \ref{intind} we can assume $f$ such that there are hyperbolic periodic points of index $k$, for every  $i\leq k\leq j$. Hence, there are $q_i, \dots, q_{j}$  hyperbolic periodic points of $f$ of indices $i, \ldots, j$, respectively. Also, by Proposition \ref{propblender} we can find a diffeomorphism $f_1$ $C^1-$close to $f$, which is not approximated by partial hyperbolic diffeomorphisms,  such that there exist blender horseshoe subsets $\Lambda_{k}$  with $ind\, \Lambda_k=k$, for every $k=i,\ldots, j-1$. By Remark \ref{csblender}, we can also assume there is a $cs-$blender horseshoe $\Lambda_{j}$ with  $ind\, \Lambda_{j}=j$.

By Proposition \ref{BC} and a result of Carballo, Moralles and Pacífico \cite{CMP}, we can also suppose that $H(q_i(f_1),f_1)=\ldots=H(q_{j}(f_1),f_1)=M$, i.e, hyperbolic periodic points of every index are dense in the whole manifold $M$. We would like to note that although the result in \cite{CMP} is in dissipative setting, provided  it is a consequence of the connecting lemma, it is still true in the volume preserving scenario. 

\begin{lemma} There exists $p\in \{ q_i(f_1),\ldots, q_{j}(f_1)\}$ and a diffeomorphism $f_2$ $C^1-$close to $f_1$ such that $f_2$ exhibits a homoclinic tangency for $p(f_2)$. \label{lemaphtg}\end{lemma}

\proof

Suppose, contrary to our claim, that every diffeomorphism in a small neighborhood of $f_1$ exhibits no homoclinic tangency for any hyperbolic periodic points $q_i(f_1),\ldots, q_{j}(f_1)$. Then, $f_1$ is in fact not approximated by diffeomorphisms exhibiting homoclinic tangency.  

Hence, by Proposition \ref{TvsD} we should have dominated splittings $TM=E_k\oplus F_k$ for $f_1$, with $dim\, E_k=k$, for every $k=i,\ldots, i+j$. Which implies we have a dominated splitting $TM=E_i\oplus 
E^1\oplus\ldots E^{j-i-1}\oplus E_j$, where $dim\, E^k=1$ for every $1\leq k\leq j-i-1$. 

Since for diffeomorphisms in $\mathcal{U}$ any hyperbolic periodic point $p$ has $ind\, p\leq j$, it follows by the same method as in Lemma 2.1 in \cite{Potrie}, that there exists $K>0$, $m\in \N$ and $0<\lambda<1$ such that every hyperbolic periodic point $p$ of a diffeomorphism $g\in \mathcal{U}$ with index $j$ and sufficiently large period one has
\begin{equation}
\prod_{l=0}^{k}\left\|\prod_{r=0}^{m-1}Dg^{-1}|E_j(g^{-lm-r}(p))\right\|\leq K\lambda^k, \text{ where } k=\left[\frac{\tau(p,g)}{m}\right].
\label{eq123}\end{equation}
To obtain this, Potrie \cite{Potrie} uses a result for uniformly contracting sequences introduced by Mañe, Lemma II.5 in \cite{Mane}.   

Now, by a known Mañe's argument also introduced in \cite{Mane},  we can find  a positive integer $n$ such that  $\|Df_1^{-n}|E_j\|$ contracts for every $x\in M$. This argument consists in use Mañe's Ergodic Closing Lemma to obtain a hyperbolic periodic point of index $j$ that doesn't satisfy equation \ref{eq123}, if $\|Df^{-n}|E_j\|$ no contracts for every $n>0$. All of these arguments, including the results for uniformly contracting sequences are done for volume preserving diffeomorphisms in \cite{AC}.

Similarly, we can prove that $\|Df_1^{n}|E_i\|$ contracts for a large enough positive integer $n$. Then $f_1$ is partial hyperbolic, which is a contradiction and then finishes the proof of lemma.   
$\hfill\square$

\begin{remark}
We point out that the above arguments give in particular a proof of Proposition \ref{CSY}.  
\label{csy}\end{remark}





Hence, let $f_2$ and $p$ given by Lemma \ref{lemaphtg}. After a perturbation, we can suppose $p$ is one of the two reference saddles of the blender-horseshoe $\Lambda=\Lambda_k(f_2)$, if $ind\, p=k$.  

By Proposition \ref{foldedmanifold}, the proof is completed by showing that:

\begin{lemma}
There is a diffeomorphism $g$ arbitrarily $C^1-$close to $f_2$ such that $W^u(p(g))$ is a folded submanifold with respect to the continuation $\Lambda(g)$  of the blender-horseshoe $\Lambda$ for $g$.
\end{lemma}

This lemma is a volume preserving version of Lemma 4.9 in \cite{BD}. 
\proof

Let $B$ denote a point of homoclinic tangency for $p(f_2)$, which we can suppose to be in $\mathbb{C}$, i.e., $B$ is in the reference cube of the blender horseshoe $\Lambda$.  As in the proof of Lemma \ref{technicallemma}, after a perturbation, we can assume $T_B W^s(p,f_2)\cap T_B W^u(p,f_2)=E^{cu}(B)$ is the one dimensional center-unstable subspace. 
Hence, let $\mathbb{V}\subset T_B W^u(p,f_2)$ be such that $T_BW^u(p,f_2)=\mathbb{V}\oplus E^{cu}$.  

Provided we have a partial hyperbolic splitting in $\mathbb{C}$,  $Df^n(B)(\mathbb{V})$ converges to $E^{uu}(p)$, if $n$ goes to infinity. Hence, if $U\subset W^u(p,f_2)$ is a small enough disk containing $B$ and  $n$ is a large enough positive integer,  
then $\mathcal{S}=f^n(U)$ is foliated by  $uu-$disks. More precisely, $\mathcal{S}=\cup_{t\in [0,1]} \mathcal{S}_t$ and $\mathcal{S}_t$ is a $uu-$disk.

We could have assumed that  $B$ is a point of quadractic homoclinic tangency.  Hence, to finish we need to analyze  two cases. 

\begin{figure}[!htb]
\vspace{0,5cm}\centering
\includegraphics[scale=1.5]{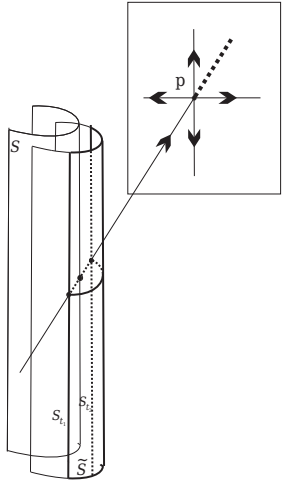}
\caption{Folded manifold: first case}\vspace{0,5cm}\label{caso a}\end{figure}

In the first one, see figure \ref{caso a},  
we can unfold the homoclinic tangency to obtain $t_1$ and $t_2$ such that $\mathcal{S}_{t_1}\cap W_{loc}^s(p,f_2)$ and $\mathcal{S}_{t_2}\cap W_{loc}^s(p,f_2)$  are non empty, and $\mathcal{S}_t$, $t_1<t<t_2$ is a $uu-$disk in between of $p$ and $q$. Therefore,  $\tilde{\mathcal{S}}=\cup_{t\in [t_1, t_2]} \mathcal{S}_t$ is a folded manifold inside the unstable manifold as we wanted.


\begin{figure}[httb]
\vspace{0,5cm}\centering
\includegraphics[scale=1.5]{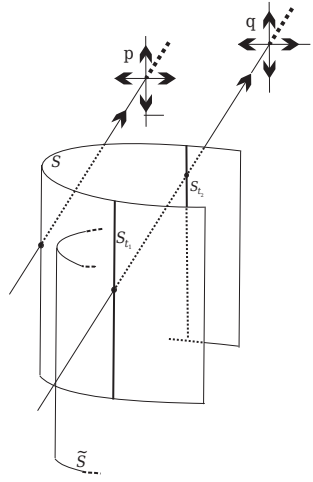}
\caption{Folded manifold: second case}\vspace{0,5cm}\label{caso b}\end{figure}

In the second case, replacing $\mathcal{S}$ by a positive iterated, should exist  $t_1$ and $t_2$ such that $\mathcal{S}_{t_1}\cap W_{loc}^s(q,f_2)$ and $\mathcal{S}_{t_2}\cap W_{loc}^s(q,f_2)$  are non empty intersections. Finally, unfolding the tangency as before, we also obtain a folded manifold. See figure \ref{caso b}.

$\hfill\square$

\section{Non existence of symbolic extensions} 

In this section we prove Propostion \ref{mainprop}, and at the end we prove Theorem \ref{simboext}. 

\vspace{0,3cm}
{\it Proof of Proposition \ref{mainprop}:}

\vspace{0,1cm}

We give the proof only for volume preserving diffeomorphisms. The general case is completely similar. 

Let $f$ be a $C^1$ volume preserving diffeomorphism and $\mathcal{U}\subset \dm$ be a neighborhood of $f$ as in the assumptions.

By Robinson \cite{Rob},  Kupka-Smale's result is still true in the volume preserving setting. Hence, there exists an open and dense subset $\mathcal{A}_{n}\subset \mathcal{U}$ of diffeomorphisms $g$ such that every periodic point of $g$ with period smaller than or equal to $n$ is hyperbolic (or elliptic if $M$ is a surface).
 Since the proposition for area preserving diffeomorphisms is contained in the Downarowicz-Newhouse's result \cite{DN},  to avoid elliptic periodic points we suppose $dim \, M\geq 3$.  

We set $\mathcal{R}_1=\cap \mathcal{A}_{n}$ which is a residual subset in $\mathcal{U}$. Let
$\mathcal{R}_{1,m,n}$ be the open set of diffeomorphisms $g$ in $\mathcal{A}_n$ where $m$ is the smallest one such that
$Per_m(g)\neq\emptyset$. Hence, 
$$
 \mathcal{A}_{n}=\cup_{j=1}^{n}\mathcal{R}_{1,j,n}.
$$

We say that an increasing sequence of finite partitions $(\alpha_n)$ over $M$ is {\it essential} for a diffeomorphism $g$  if
\begin{itemize}
\item [1.] $diam(\alpha_k)\rightarrow 0$ when $k\rightarrow \infty$, and
\item [2.] $\mu(\partial\alpha_k)=0$ for every $\mu\in\mathcal{M}(g)$. Here, $\partial\alpha_k$
denotes the union of boun\-da\-ries of all elements of the partition  $\alpha_k$.
\end{itemize}

By Proposition 4.1 in \cite{DN}, we can assume there exists an increase sequence of partitions  $(\alpha_n)$ over $M$ and a residual subset $\mathcal{R}_2\subset \mathcal{U}$ such that $(\alpha_n)$ is an essential sequence for every diffeomorphism $g\in \mathcal{R}_2$. From now on we fix this partition. 

Recall, a $f-$invariant periodic set $\Lambda$ with basis $\Lambda_1$  is {\it subordinated} to a finite partition $\alpha$ if for each  positive integer $n$, there exists an element $B_{n}\in\alpha$ such that
$f^{n}(\Lambda_1)\subset B_{n}$. In particular, if $\Lambda$ is subordinated to $\alpha$ and $\mu\in \mathcal{M}(f|\Lambda)$ then $h_{\mu}(\alpha)=0$.



 Let $\Lambda$ be the hyperbolic basic set which exhibits robust homoclinic tangency for $f$. For every $g\in \mathcal{U}$, we denote by 
$H(\Lambda(g))$ the set of hyperbolic periodic points of $g$ homoclinically related with $\Lambda(g)$, and  we set
$$
H_n(\Lambda(g))=H(\Lambda(g))\cap Per_n(g).
$$

If $p$ is a hyperbolic periodic point of $g$, $|\mu(p,g)|<1<|\lambda(p,g)|$ denote the two eigenvalues of
$Dg^{\tau(p,g)}$ nearest one, i.e., if $\nu$ is an eigenvalue of $Dg^{\tau(p,g)}$ then either $|\nu|\leq |\mu(p,g)|$ or $|\nu|\geq |\lambda(p,g)|$. 
We define  
$$
\chi(p,g)=\frac{1}{\tau(p,g)}\log\min\{|\lambda(p,g)|,\,|\mu(p,g)|^{-1}\}, $$ for every hyperbolic periodic
point $p$ of $g$.

Finally, for any positive integer $n$,  we say that a diffeomorphism $g\in \mathcal{U}$ satisfies  pro\-per\-ty
$\mathcal{S}_n$ if for every $p\in H_n(\Lambda(g))$  {\it
\begin{itemize}
\item [1.] There exists a hyperbolic basic set of zero dimension $\Lambda(p,n)$ for $g$ such that
\begin{equation*}
\Lambda(p,n)\cap\partial\alpha_n=\emptyset \; \text{ and } \; \Lambda(p,n)\, \text{ is subordinate to }
\,\alpha_n.
 \label{cond 1}\end{equation*}

\item[3.] There exists an ergodic measure $\mu\in \mathcal{M}(f|\Lambda(p,n))$ such that
\begin{equation*}
h_{\mu}(g)> \chi(p,g)-\frac{1}{n}. \label{cond 3}\end{equation*}

\item[4.] For every ergodic measure  $\mu\in\mathcal{M}(f|\Lambda(p,n))$, we have
\begin{equation*}
\rho(\mu,\mu_p)<\frac{1}{n},
 \label{cond 4}\end{equation*}
where $\mu_p$ is the dirac measure on the orbit of $p$.

\item[5.] For every periodic point $q\in \Lambda(p,n)$, we have
\begin{equation*}
\chi(q,g)>\chi(p,g)-\frac{1}{n}. \label{cond 5}\end{equation*}
\end{itemize}}

For positive integers  $m\leq n$, let $\mathcal{D}_{m,n}\subset \mathcal{R}_{1,m,n}$ be the subset of
diffeomorphisms  satisfying property $\mathcal{S}_n$.

\vspace{0,3cm}
{\bf Claim:}  $D_{m,n}$ is open and dense in $\mathcal{R}_{1,m,n}$.

\vspace{0,2cm}

This Claim is a conservative version of Lemma 3.3 in \cite{CT}, which is an extension for symplectic diffeomorphisms of Lemma 5.1 in \cite{DN}.  
By Claim the proof is similar in spirit to Theorem 1.3 in \cite{DN}. See also \cite{CT}. However, just for sake of completeness we will give a sketch here. 

We remark first if $(\alpha_n)$ is an essential sequence of partitions for $f$, then  for any $k$ fixed we set
$$h_k(\mu)=h_{\mu}(\alpha_k),$$ which is an infimum of continuous functions over $\mathcal{M}(f)$. 

The following
proposition relates the entropy structure of a diffeomorphism and non existence of symbolic extensions.  It was also proved in
\cite{DN}.


\begin{lemma} Let $f\in \mathcal{R}_2$.  
If
there exists a compact subset $\mathcal{E}\subset \mathcal{M}(f)$ and $\rho>0$ such that
$$
\limsup_{\nu\rightarrow \mu, \nu\in \mathcal{E}} h_{\nu}(f)-h_k(\nu)\geq \rho, \text{ for every } \mu\in
\mathcal{E} \text{ and } k\geq 0,
$$
then $f$ has no symbolic extension. 
\label{lema 1}\end{lemma}

We set $ \mathcal{R}=\cap_{n\geq0}\cup_{m=0}^n D_{m,n}\cap \mathcal{R}_2$, which is a residual subset in
$\mathcal{U}$ by  Claim.


Now, let $f\in\mathcal{R} $ and define  $\chi(f)=\sup\{\chi(p,f), p\in H(\Lambda(f))\}$.  We denote by
$\mathcal{E}$  the closure of
$$
\mathcal{E}_1=\{\mu_p;\, p\in H(\Lambda(f)) \text{ and } \chi(p,f)>\chi(f)/2\}.
$$


Since $f\in \mathcal{R}$, for any periodic point $p$ such that $\mu_p\in
\mathcal{E}_1$, it follows there exist ergodic measures $\nu_n\rightarrow \mu_p$ such that $h_{\nu_n}(f)>\chi(f)/2$.
Moreover, 
 $\nu_n\in \mathcal{M}(f|\Lambda(p,n))$, 
by Sigmund \cite{SI}, $\nu_n$ is approximated by hyperbolic
periodic measures also supported in the hyperbolic set $\Lambda(p,n)$, which by item 4 of property $S_n$ should
be in $\mathcal{E}_1$. Hence, $\nu_n\in \mathcal{E}$ for every $n$.

Therefore, defining $\rho=\chi(f)/2>0$ and $\mathcal{E}$ as before,  by Lemma \ref{lema 1} $f$ has no symbolic extensions. 

\vspace{0,2cm} {\it Proof of Claim:}  

As in the proof of Technical Proposition in \cite{CT}, the procedure is to find Newhouse's snakes (see Remark \ref{snake} for a definition) after a perturbation of a diffeomorphism exhibiting a homoclinic tangency,  to obtain from them nice hyperbolic sets satisfying the conditions in property $S_n$. However, to find a diffeomorphism satisfying property $S_n$ it's necessary to have an argument to obtain Newhouse's snakes related to any arbitrary hyperbolic periodic point. To do this we will use robustness of homoclinic tangency.

We denote by $\Lambda$ the hyperbolic set of $f$ exhibiting robust homoclinic tangency.

Let $g\in \mathcal{A}_{n}$ be an arbitrary diffeomorphism. By definition of $\mathcal{A}_n$, there exists a
small neighborhood $\mathcal{V}$ of $g$ where the cardinality of periodic points of period smaller than or equal to $n$ is constant. 

Let $p\in H_n(\Lambda(g))$. Since $\Lambda(g)$ has a robust homoclinic tangency and $p$ is homoclinic related with $\Lambda(g)$, after a perturbation we can
assume that $g$ exhibits a homoclinic tangency for the hyperbolic periodic point $p$. By regularization result of Ávila \cite{Avila} we can suppose that $g$ is $C^{\infty}$ and then using pasting lemma we can assume that  $g^{\tau(p,g)}=Dg^{\tau(p,g)}$  in some neighborhood of $p$ (in local coordinate), say $V$. That is,  $g^{\tau(p,g)}|(V\cap g^{-\tau(p,g)}(V))$ is linear. Hence, $W^s_{loc}(p,g)$ and $W^u_{loc}(p,g)$ coincide with stable and unstable directions restrict to $V$.

For simplicity we suppose $p$ is a hyperbolic fixed point of $g$. 

Let $q$ be the point of homoclinic tangency between $W^s_{loc}(p,g)$ and $W^u(p,g)$, such that $q\in V$ and $g^{-1}(q)\not\in V$.  Hence,
 there is a small neighborhood $U\subset V$ of $q$ such that $g^{-1}(U)\cap V=\varnothing$. We denote by $D$ the connected component of
$W^u(p,g)\cap U$ that contains $q$. For convenience we suppose $dim(TqW^s_{loc}(p,g)\cap T_qW^u(p,g))=1$, which can be done after a perturbation.

We look to $U$ in a local coordinate,  being $q$ the zero of this chart, and we consider the following splitting of space $\R^d=T_qD\oplus T_qD^{\bot}$. Since $D\cap U$ is an open disc inside $W^u(p,g)$, it follows that $D$ is a graph of a map $r:T_{q}D\rightarrow T_{q}D^{\bot}$, which  is so regular as $g$. That is,  $D=(x,r(x))$ is the graphic of a $C^{\infty}$ map $r$. Moreover, such map $r$ is such that $Dr(q)$ is zero. Defining $\phi(x,y)=(x,y-r(x))$, if $U$ is small enough, this map is is a $C^{\infty}$ volume preserving diffeomorphism from $U$ to $\phi(U)$, and is $C^1$ close to identity in a small enough neighborhood of $q$.
Therefore, by pasting lemma we can find a $C^{\infty}$ volume preserving diffeomorphism $h$ on $M$, $C^1$-close to identity such that $h=\phi$ in some small neighborhood of $q$, and $h=Id$ outside $U$.

We define $g_1=h\circ g$ which is a $C^1$ perturbation of $g$.  Note, $g_1$ is such that  $T_qD\cap U\subset W^u(p,g_1)$. Since $g^{-1}(U)\cap V=\emptyset$,  it follows that $g_1=g$ in $V$ and so $g_1|V$ is still linear, which implies  $W^s_{loc}(p,g_1)\cap U=W^s_{loc}(p,g)\cap U=E^s(p,g)\cap U$. Hence, provided $q$ was a non-transversal homoclinic point of $p$ for $g$ (i.e., $T_q D\cap E^s(p,g)$ is non trivial),  $g_1$ exhibits an interval of homoclinic tangency containing $q$.


Let $I$ be this interval of homoclinic tangency. By a coordinate change in $U$,  we can suppose that $W^s_{loc}(p,g_1)\cap U=E^s(p,g_1)\cap U\subset \R^s\times \{0\}^u$, and $I\subset\{(x_1,0,...,0), \, -3a\leq x_1\leq 3a\}$, for some $a>0$ small enough. We are now considering the usual coordinates in
$\R^{d}$.

Let $N$ be any large positive integer. By the same method as in the construction of $h$, for any $\delta>0$ small enough we can find a
volume preserving diffeomorphism $\Theta:M\rightarrow M$, $\delta-C^1$ near $Id$, such that $\Theta=Id$ in the complement of $B(q,2a)$ and
$$
\Theta(x,y)=\left(x_1,...,x_s,\,y_1+A\cos\frac{\pi x_1 N}{2a},\, y_2,...,y_u\right), \text{ for } (x,y)\in B(0,a)\subset U,
$$
for $A=\displaystyle\frac{2Ka\delta }{\pi N}$, where $K$ is a constant depending only on the local coordinate on  $U$. 

We define $g_2=\Theta\circ g_1$ which is $\delta-C^1$ close to $g_1$ and moreover $g_2=g_1$ in the complement of $g_1^{-1}(U)$.  Note, $g_2$ exhibits $N$ transversal homoclinic points for $p$ inside $U$. More precisely, these points belong to $g_2(g^{-1}(I))$. 

\begin{remark}
This kind of perturbation is the so called Newhouse's snake.  
\label{snake}\end{remark}

We remark that $g_2$ depends on $N$, but by abuse of notation, we use the same letter $g_2$ for every $N$.


From now on we will look to $V$ in local coordinate. Moreover, we assume  $E_{p}^s=\R^{s}\times \{0\}^u$ and
$E_{p}^u=\{0\}^s\times\R^{u}$, where $E_{p}^s$ and $E_{p}^u$ are the stable and unstable directions of $p$ with dimensions $s$ and $u$, respectively. Observe $g_2$ is linear in $V$ since $g_2$ is equal to $g_1$ in $V$. 

For any  positive large integer  $t$, we can define a small rectangle $D_t=D^s\times D^u_t$, where $D^s=W^s_{loc}(p,g_2)\cap U$ and $D^u_t$ is a small disk in $\{(0,\ldots,0,y_1,\ldots,y_n),\, y_i\in\R^+ \text{ and } |y_i|<A/4\}$, such that $g_2^t(D_t)$ is a disk $A/4-C^1$ close to the connected component of $W^u(p,g_2)\cap U$ containing the $N$ transversal homoclinic points of $g_2$ in $g_2(g^{-1}(I))$. We fixe $t$ as being the smallest possible one such that $D_t$ is defined.  It is clearly that $t$ depends on $N$, and goes to infinity if $N$ goes. 

Since $N$ is large, it follows that $A$ is so small which implies $D_t$ is such that $g_2(D_t)\cap D_t$ has $N$ disjoint connected components. Moreover, taking $N$ larger if necessary, the maximal invariant set in $D_t$ for $g_2^t$
$$
\tilde{\Lambda}(p,N)=\bigcap_{j\in\Z} g_2^{tj}(D_t)
$$
is a hyperbolic set. 

Let $$\Lambda(p,N)=\bigcup_{0\leq j\leq t} g_2(\tilde{\Lambda}((p,N)))$$ be the hyperbolic periodic set of $g_2$ induced by $\tilde{\Lambda}(p,N)$.

Now, if $n$ is an arbitrary large positive integer $n$, we can proceed in the same way as in the proof of technical proposition in \cite{CT} to find a large positive integer $N$, such that  $\Lambda(p,N)(\tilde{g})$ satisfies all items of property $S_n$ for every diffeomorphism $\tilde{g}$ $C^1-$close to $g_2$.


Finally, since 
the cardinality of $Per_n$ is finite and constant for diffeomorphisms in $\mathcal{V}$, we can repeat the same process finitely many times to obtain an open set $C^1$-close to $g$ of diffeomorphisms satisfying property $S_n$. 

$\hfill\square$

We finish the paper with the proof of Theorem \ref{simboext}.

\vspace{0,2cm}
{\it Proof of Theorem \ref{simboext}:}

Since $\dm$ is a separable space, it follows there exists an enumerable dense subset $\{f_1,\ldots,f_n,\ldots\}$ of diffeomorphisms in $\dm$. 

If  $f_i$ is not partial hyperbolic, then by Theorem \ref{simboext} we can suppose $f_i$ exhibits a robust homoclinic tangency, after a perturbation. Then, by Proposition \ref{mainprop} there exist a neighborhood $\mathcal{U}_i$ of $f_i$ and a residual subset $\mathcal{R}_i\subset \mathcal{U}_i$ such that every diffeomorphism $g\in \mathcal{R}_i$ has no symbolic extensions. 
Hence, $\mathcal{R}_i$ contains an enumerable intersection of open and dense subset in $\mathcal{U}$, say $\mathcal{R}_i=\cap \tilde{\mathcal{B}}_n^i$.
We define $\mathcal{B}_n^i=\tilde{\mathcal{B}}_n^i\cup (cl(\mathcal{U}_i))^c$, which is in fact an open and dense subset of $\dm$.  

If $f_i$ is partial hyperbolic we define $\mathcal{B}_n^i$ as being the set of partial hyperbolic diffeomorphisms and the ones that are not approximated by partial hyperbolic diffeomorphisms, which is also an open and dense subset of $\dm$. 

Hence, we define 
$$
\mathcal{R}=\bigcap_{i,n\in \N} \mathcal{B}_n^i,
$$  
which is a residual subset in $\dm$. 
Finally, note by construction that if $f\in \mathcal{R}$ is not partial hyperbolic, then $f$ has no symbolic extensions. Which proves the theorem. 

$\hfill\square$

\vspace{0,3cm}
{\bf Acknowledgements:} I would like to thanks S. Crovisier and E. Pujals for useful conversations, in particular   the difficulties to obtain results in this direction for higher topologies. To A. Arbieto, V. Horita and A. Tahzibi for their mathematical support and constant encouragement. To IBILCE-UNESP, IMPA  and The Abdus Salam ICTP for the kind hospitality during preparation of this work. To FAPESP-Brazil for its financial support.

\bibliographystyle{amsplain}

\vspace{0,3cm} \flushleft
Thiago Catalan\\
Faculdade de Matemática \\
Universidade Federal de Uberlândia \\
34-32309442 Uberlândia-MG, Brazil\\
E-mail: tcatalan@famat.ufu.br

\end{document}